\documentclass[11pt,a4paper]{article}
\usepackage[cp1251]{inputenc}
\usepackage[english,russian]{babel}
\usepackage{amsmath, nccmath}
\usepackage{amssymb}
\usepackage{euscript}
\usepackage{wrapfig}
\usepackage{cite}
\usepackage{amscd}
\usepackage{amsthm}
\usepackage{color}
\usepackage{graphics}
\usepackage{pdflscape}
\usepackage{epsfig}
\usepackage{mathrsfs}
\usepackage{latexsym}
\usepackage{srcltx}

\usepackage{textcomp}
\AtBeginDocument{}

\makeatletter
\renewcommand\@seccntformat[1]{\csname the#1\endcsname\hspace{0.5em}}
\makeatother

\newtheorem{theorem}{Теорема}[section]
\newtheorem*{theorem*}{Теорема}
\newtheorem{lemma}[theorem]{Лемма}
\newtheorem{prop}[theorem]{Предложение}

\newtheorem{remark}[theorem]{Замечание}
\newtheorem{quation}[theorem]{Вопрос}

\begin{document}
\begin{center}
{\Large \bfseries Feuerbach's and Poncelet's theorems meet in space}\par
{\it (On the occasion of their bicentennial)}

\medskip
{E.\,A. Avksentyev}

\medskip
{\it December 29, 2022}
\end{center}
\begin{abstract}
We propose 3D generalizations of the Feuerbach theorem: the first one deals with a tetrahedron analogue of the Euler circle, the second one is done by means of an {\guillemotleft}up-in-ex-touch{\guillemotright}  construction.  Then we give a geometric proof of the Grace theorem (a classical, but still not well-known,  3D Feuerbach theorem) and show its relation to the Poncelet closure theorem.  Our elementary proof  is based on the properties  of imaginary generators on a sphere and of isotropic tangents to a conic. Then we show that the Grace theorem implies the  Laguerre theorem, which generalises the Euler-Chapple formulas. Further, we consider a spatial analog of Poncelet’s theorem. It is also proved that the Grace spheres touch a fixed sphere under the Poncelet rotation of bicentric tetrahedron. Finally, the lifting to the three-dimensional space  provides a new proof of Feuerbach’s theorem and, perhaps, the shortest proof of Euler-Chapple formulas.

\end{abstract}
%


\section*{Введение}

Данная работа посвящена двум знаменитым геометрическим теоремам, кажется никак не связанным между собой, разве что они были опубликованны в один год двести лет назад~\cite{Poncelet,Feuerbach}. Приведем их формулировки

\begin{theorem*} [Feuerbach, 1822]
\label{Feuerbach}
Окружность девяти точек произвольного треугольника касается его вписанной и трех вневписанных окружностей.
\end{theorem*}

\begin{theorem*} [Poncelet, 1822]
\label{th_Poncelet}
Пусть для двух данных коник существует вписано-описанный в них многоугольник. Тогда этот многоугольник может динамически «вращаться» около данных коник, оставаясь вписано-описанным в них.
\end{theorem*}


У обеих теорем есть масса обобщений, но пространственные аналоги, насколько нам известно, имеются только у теоремы Понселе. Их довольно много (см., например, ~\cite{Fontene,Grace,Griffiths,Protasov}) и среди них есть множество замечательных, но малоизвестных результатов.

Задача трехмерного обобщения теоремы Фейербаха поставлена еще более ста лет назад в монографии Кулиджа~\cite{C}:
\begin{center}
\begin{minipage}[h]{0.85\linewidth}
«{\it The geometry of the tetrahedron lags far behind that of the triangle... Is there an analogue to Feuerbach's~theorem? Above all what corresponds to the
Hart systems? ...These difficult but important and interesting questions offer ample~scope~for~serious~work}»~(p.\,247).
\end{minipage}
\end{center}

Теорема Фейербаха содержит в себе два удивительных геометрических факта. Первый состоит в том, что четыре замечательные окружности треугольника~-- вписанная и три вневписанные~-- имеют общую касательную окружность. Второй же заключается в том, что эта общая касательная окружность является еще и окружностью девяти точек, которая и без того сама по себе замечательна.

Первая попытка найти аналог теоремы Фейербаха в пространстве приводит к вопросу: {\it существует ли сфера, которая касалась бы вписанной и вневписанных сфер?}

\begin{figure}[htb]
\center
\includegraphics[scale=0.4]{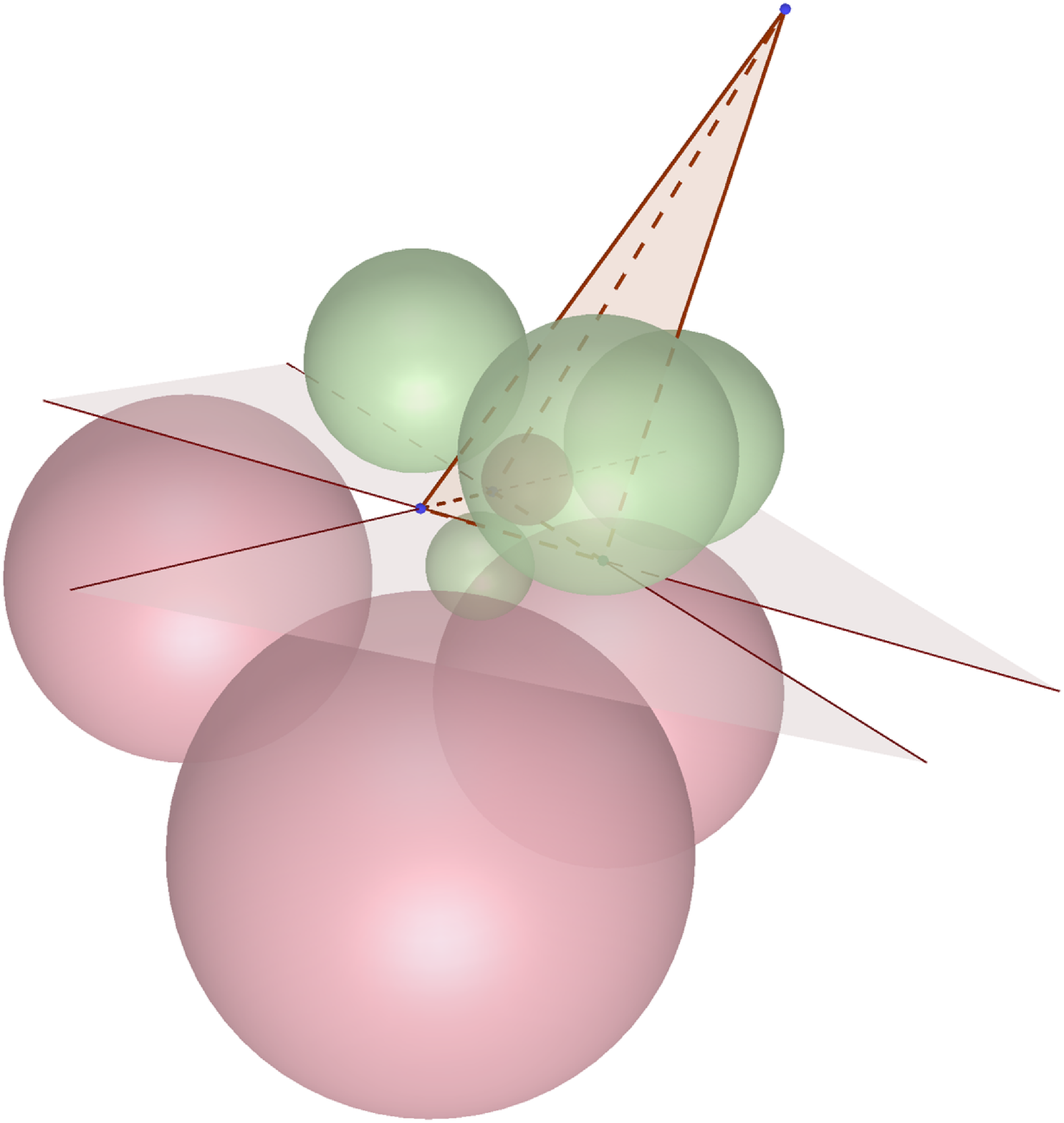}
\caption{\small{\it Восемь касательных сфер тетраэдра}}
\label{fig_In-ex-spheres}
\end{figure}
Но здесь нас ожидает первый «сюрприз»: у произвольного тетраэдра кроме обычных четырех вневписанных сфер, аналогичных трем вневписанным сферам треугольника, существует еще три {\it дважды-вневписанные} сферы или {\it чердачные} (от англ. «roof»), как они названы в~\cite{Berger} (см. также~\cite{Zaslavsky}).

Т.е., всего существует целых восемь сфер~(см.~рис.~\ref{fig_In-ex-spheres}), касающихся граней тетраэдра!
Назовем их {\it касательными сферами}.
Было бы слишком оптимистично ожидать, что все восемь касательных сфер могли бы касаться одной сферы. И действительно, ответ на поставленный вопрос оказывается отрицательным: {\it в общем случае произвольного тетраэдра такой сферы не существует.}

Проверить это очень легко: для этого достаточно рассмотреть лишь один подходящий тетраэдр, например, равногранный. И нет сомнений, что такой знаток геометрии как Кулидж хорошо знал, что такой сферы в общем случае нет. Однако, он все-таки поставил вопрос поиска трехмерных аналогов теоремы Фейербаха, находя его важным, интересным и открывающим «широкие возможности для серьезной работы».

В каком же направлении искать тогда аналоги теоремы Фейербаха в пространстве? Кажется, что осталась лишь задача описания частных случаев тетраэдров, у которых существует сфера, касающаяся внутренним или внешним образом пяти, шести, семи или всех восьми касательных сфер. В~работе~\cite{Lewis} есть некоторое продвижение в этой задаче и для существования такой сферы получены аналитические условия в специальных связанных с тетраэдром пентасферических координатах. К сожалению, эти условия весьма громоздкие и из них совершенно не ясно, существуют ли такие тетраэдры и как они устроены. Таким образом, задача в такой постановке остается незакрытой.

Возникает еще идея поискать пространственный аналог теоремы Фейербаха в таком направлении: {\it существует ли окружность, действительная или мнимая, которая касалась бы всех восьми касательных сфер?} Кажется маловероятным, что ответ мог бы быть положительным, но задача представляется интересной.

Оставив пока эти вопросы, мы приведем далее целых три трехмерных аналога теоремы Фейербаха.

Первый аналог, который мы хотим предложить в \S~\ref{ch_1analog_cones} в качестве трехмерного обобщения теоремы Фейербаха, является довольно интересным фактом. У него очень простое доказательство, которое, тем не менее, раскрывает связь этой конструкции с неевклидовой геометрией и приводит к трехмерному обобщению окружности Эйлера. Поэтому из трех аналогов этот наиболее аутентичен.

Второй является очень красивой теоремой геометрии тетраэдра, открытой сто двадцать пять лет назад, но, кажется, до сих пор малоизвестной. Ее единственное оригинальное доказательство столь сложно, что есть целая статья с его реконструкцией. В \S~\ref{ch_Grace} мы получим элементарное доказательство этой теоремы, в котором обнаружится ее связь с теоремой Понселе. Второй аналог выглядит наименее аутентичным, но на наш взгляд, он ближе и роднее к теореме Фейербаха, чем другие два.

Третий аналог представляет из себя довольно интересную конструкцию касающихся сфер, которую мы назвали «up-in-ex-touch»-конструкция. Мы приведем ее в конце \S~\ref{ch_up-in-ex-touch}, в котором мы также получим, возможно, самое короткое доказательство формул Эйлера-Чаппла.

С помощью теоремы Грейса мы в \S\ref{ch_Laguerre} получим короткое и простое доказательство теоремы Лагерра, обобщающей формулы Эйлера-Чаппла. \S\ref{ch_R3_Euler-Chapple} посвящен трехмерному аналогу формул Эйлера-Чаппла.

Далее в \S\ref{ch_R3_Ponc} мы рассмотрим пространственные аналоги теоремы Понселе. Мы покажем, что при вращении Понселе вписано-вневписанного тетраэдра его сферы Грейса касаются некоторой фиксированной сферы.

В конце, совершая «выход в пространство», мы дадим новое доказательство теоремы Фейербаха.

\section{Первый аналог теоремы Фейербаха для тетраэдра}
\label{ch_1analog_cones}
Итак, рассмотрим произвольный тетраэдр общего вида, у которого имеется восемь касательных сфер. В качестве первого аналога теоремы Фейербаха для тетраэдра предлагаем следующую теорему.

\begin{theorem}
\label{R3_Feuerbach}
Существует четыре круговых конуса, каждый из которых касается всех восьми его касательных сфер.
\end{theorem}

{\tt Доказательство.}
Рассмотрим сферу $\zeta_D$ с центром в вершине $D$ тетраэдра $ABCD$ и спроектируем из центра $D$ на сферу $\zeta_D$ все восемь касательных сфер. Их проекциями будут четыре окружности на сфере $\zeta_D$, поскольку каждая пара гомотетичных относительно $D$ сфер спроектируются в одну и ту же окружность. Эти четыре окружности касаются сторон сферического треугольника, стороны которого являются проекциями плоскостей трехгранного угла при вершине $D$. По теореме Фейербаха для сферического треугольника существует окружность, касающаяся этих четырех окружностей. Конус с вершиной $D$, содержащий эту окружность, очевидно удовлетворяет утверждению теоремы. Такой конус есть у каждой вершины.
{\hfill $\Box$}
\medskip

Теорема Фейербаха в сферической геометрии, в той облегченной форме, которую мы использовали в доказательстве, равносильна теореме Харта~(см.~\cite{C}). Таким образом, в какой-то степени мы ответили на оба вопроса Кулиджа, которые мы цитировали во введении. На самом деле, можно продвинуться еще дальше в этом направлении, если применить результат Акопяна~\cite{Akopyan}, в котором он нашел такие свойства окружности Харта, которые во многом аналогичны свойствам окружности девяти точек. Хотя в~\cite{Akopyan} все утверждения формулируются для плоскости Лобачевского, но мы их естественным образом адаптируем применительно к трехгранным углам нашего тетраэдра.

{\it Избытком} трехгранного угла называется величина, равная разнице между суммой его двухгранных углов и $180^{\circ}$. {\it Медиатором} трехгранного угла назовем плоскость, содержащую его ребро и делящую его на два трехгранных угла с равными избытками. При рассмотренной выше проекции трехгранного угла на сферу медиатор переходит в сферическую чевиану, делящую пополам площадь соответственного треугольника (в~\cite{Akopyan} эта чевиана называется биссектором или биссекторным отрезком). Три медиатора пересекаются по прямой, которую можно назвать {\it псевдоцентроидалью}, поскольку ей соответствует псевдоцентроид сферического треугольника.

Четыре прямые из одного пучка назовем {\it вписанной четверкой}, если все они являются образующими одного кругового конуса. Следующее утверждение является аналогом Леммы~5 из~\cite{Akopyan}.

\begin{prop}
\label{th_pseudo_altitudes}
Пусть $a, b, c$~-- ребра трехгранного угла с вершиной $D$. Тогда существует единственная тройка прямых $\mathfrak{h}_a, \mathfrak{h}_b, \mathfrak{h}_c$, лежащих в плоскостях $\langle bc\rangle,\,\langle ac\rangle, \langle ab\rangle$ соответственно, таких что четверки $\{a, b, \mathfrak{h}_a, \mathfrak{h}_b\}$;\, $\{a, c, \mathfrak{h}_a, \mathfrak{h}_c\}$;\, $\{b, c, \mathfrak{h}_b, \mathfrak{h}_c\}$ являются вписанными.
\end{prop}

Плоскости $a\mathfrak{h}_a,\,b\mathfrak{h}_b,\,c\mathfrak{h}_c$ являются аналогами так называемых {\it псевдовысот}, которым в~\cite{Akopyan} дается еще и другое определение через углы. Эти три плоскости пересекаются по общей прямой, назовем ее {\it псевдоортоцентралью} по аналогии с псевдоортоцентрами гиперболических треугольников.

Круговой конус, содержащий все три ребра трехгранного угла в качестве своих образующих, назовем {\it описанным}.

В~\cite[\S\S\,4-6]{Akopyan} показано, что основания трех псевдовысот и трех биссекторных чевиан лежат на одной окружности. Центр этой окружности лежит на одной прямой с центром описанной, псевдоцентроидом и всевдоортоцентром. Сформулируем аналогичные утверждение для тетраэдра.

\begin{theorem}[Конус Эйлера трехгранного угла]
\label{th_Pseudo_Euler}
У любого трехгранного угла основания трех его медиаторов и трех его псевдовысот лежат на одном круговом конусе.
\end{theorem}

\begin{theorem}[Плоскость Эйлера трехгранного угла]
\label{th_Pseudo_Euler}
У произвольного трехгранного угла четыре прямые~-- псевдоцентроидаль, псевдоортоцентраль, ось описанного конуса и ось конуса Эйлера~-- лежат в одной плоскости.
\end{theorem}

Главным же результатом работы~\cite{Akopyan} является гиперболический аналог теоремы Фейербаха, согласно которому окружность Эйлера гиперболического треугольника касается его вписанной и трех вневписанных окружностей. Применительно к тетраэдру мы получаем следующее усиление Теоремы~\ref{R3_Feuerbach}

\begin{theorem}[Аналог теоремы Фейербаха для тетраэдра]
\label{R3_Feuerbach'}
Четыре конуса Эйлера трехгранных углов тетраэдра касаются всех восьми его касательных сфер.
\end{theorem}

Отметим несколько вопросов, которые возникают в связи с рассмотренными конструкциями.

\begin{quation}
\label{quat_incident}
Инцидентны ли какие либо из следующих четверок замечательных прямых тетраэдра: псевдоцентроидали, псевдоортоцентрали, оси четырех описанных конусов, оси четырех конусов Эйлера?
\end{quation}

\begin{quation}
\label{quat_quadric} Существуют ли еще какие-либо квадрики, касающиеся всех касательных сфер, отличные от четырех конусов Эйлера и четырех плоскостей граней?
\end{quation}

\begin{quation}
\label{quat_biquadratiques} Любые три конуса общего положения пересекаются в восьми точках. Не окажется ли так, что четыре конуса Эйлера тетраэдра имеют восемь общих точек? Есть ли какие-то примечательные свойства биквадратических кривых, по которым пересекаются конусы Эйлера?
\end{quation}

\section{Теорема Грейса как трехмерный аналог теоремы Фейербаха}
\label{ch_Grace}

Более ста лет назад, британский математик Джон Хилтон Грейс в своей работе~\cite{Grace1897} открыл и доказал следующее замечательное свойство касательных сфер тетраэдра.
\begin{theorem} [Grace, 1897]
\label{th_Grace}
Касательные сферы тетраэдра $ABCD$ могут быть разбиты на четыре пары так, что парные сферы гомотетичны с центром $D$, и для каждой пары существует касающаяся их сфера, проходящая через вершины $A, B, C$.
\end{theorem}

\begin{remark}
\label{rem_Grace1}
Все касательные сферы можно разбить на две группы по четыре сферы. В одну входят вписанная и три дважды-вневписанные сферы, а в другую~-- четыре вневписанные. Любые две сферы из разных групп гомотетичны относительно одной из вершин тетраэдра. Для каждой такой пары сфер существует единственная касающаяся их сфера Грейса, которая проходит через вершины грани, противоположной вершине, относительно которой данная пара касательных сфер гомотетична. Таким образом, всего получается шестнадцать сфер Грейса: для каждой из четырех граней тетраэдра через ее вершины проходит четыре различные сферы Грейса.
\end{remark}

Теорема Грейса связывает касательные сферы тетраэдра с замечательными точками, его вершинами, с помощью общих касающихся их сфер. Это ее сближает с теоремой Фейербаха, с которой она, на наш взгляд, сравнима по красоте и имеет некоторое сходство. В этом смысле, можно было бы считать теорему Грейса трехмерным аналогом теоремы Фейербаха.


\begin{figure}[htb]
\center
\includegraphics[scale=1]{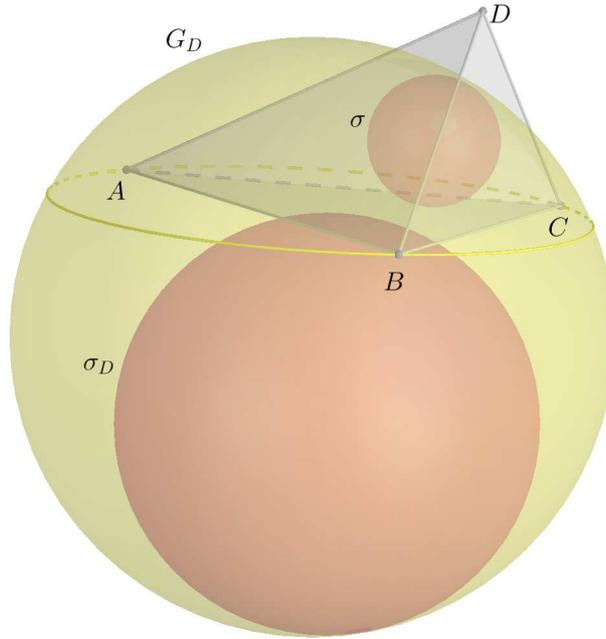}
\caption{\small{\it Сфера Грейса $G_D$, касающаяся вписанной сферы $\sigma$, вневписанной сферы $\sigma_D$ и проходящая через $A, B, C$.}}
\label{fig_Grace_sphere}
\end{figure}

В недавней статье~\cite{Maehara-Martini} Maehara и Martini замечают, что «по-видимому, эта теорема малоизвестна и до сих пор не имеет элементарного доказательства». В качестве результата они приводят такое доказательство, но лишь для частного случая триортогонального тетраэдра, пользуясь при этом аналитической техникой.

Оригинальное же доказательство Грейса очень красивое и геометрическое, но довольно трудное. Поскольку Грейс дал лишь его набросок, Maehara и Tokushige в работе~\cite{Maehara-Tokushige} подробно реконструировали это доказательство.

Мы получим элементарное и достаточно короткое геометрическое доказательство теоремы Грейса, но сначала напомним некоторые определения и факты проективной геометрии.

Пусть $\mathbb{E}^3$~-- вещественное трехмерное евклидово пространство. Мы будет рассматривать его проективное пополнение «бесконечно удаленной» плоскостью. Эта модель проективного пространства получается переходом от декартовых координат $(x, y, z)$ в $\mathbb{E}^3$ к однородным координатам $(x: y: z: w)$, в которых бесконечно удаленной плоскости соответствуют точки с координатами $(x: y: z: 0)$. Кроме того рассмотрим комплексификацию пространства, позволяя координатам принимать комплексные значения. Добавленные точки будем называть {\it мнимыми}.

Записывая в однородных координатах $(x: y: z: w)$ общее уравнение сферы
$$x^2 + y^2 + z^2 + 2axw + 2byw + 2czw + dw^2 = 0,$$
легко видеть, что она пересекает бесконечно-удаленную плоскость $w=0$ по кривой
$$x^2 + y^2 + z^2 = 0,\ w=0,$$
которая является общей для всех сфер. Она называется {\it абсолютной окружностью.}

Всякая плоскость пересекает абсолютную окружность в двух сточках~-- {\it {круговых точках}} этой плоскости. В однородных координатах $(x: y: z)$ на плоскости ее круговыми точками являются точки $I = (1: i: 0)$ и $J = (1: -i: 0)$. Все окружности плоскости проходят через ее круговые точки и каждая коника плоскости, проходящая через ее круговые точки, является окружностью~(см. \cite[\S\,4$\cdot$8]{Sommerville}).

Прямая, пересекающая абсолютную окружность, называется {\it изотропной}. Каждая такая прямая является либо бесконечно удаленной, либо мнимой.

\begin{prop}[{\cite[Гл.\,12, \S\,2]{Finikov}}]
\label{prop_isotrop}
Касательные к невырожденной конике, проведенные из любого ее фокуса, являются изотропными.
\end{prop}

Таким образом, каждая прямая, проходящая через фокус коники и круговую точку ее плоскости, является изотропной. Для окружности это означает, что касательные из ее центра проходят через круговые точки.

{\it Образующей}
 квадрики называется прямая, которая целиком принадлежит поверхности этой квадрики. В комплексном проективном пространстве все невырожденные квадрики эквивалентны.

\begin{prop}[{\cite[\S \,2]{Griffiths}}]\
\label{prop_generating_lines}
\begin{description}
  \item[(i)] Через каждую точку невырожденной квадрики проходят ровно две образующие, действительны или мнимые. Касательная плоскость пересекает квадрику по двум образующим, проходящим через точку касания.
  \item[(ii)] Все образующие квадрики распадаются на два семейства таким образом, что любые две образующие из одного семейства не пересекаются, а любые две образующие из разных семейств пересекаются. Через любую точку образующей одного семейства проходит единственная образующая другого семейства.
  \item[(iii)]\label{iii} Любая плоскость, проходящая через образующую квадрики, касается этой квадрики в некоторой точке этой образующей.
\end{description}
\end{prop}

Пусть даны две сферы $\gamma$ и $\eta$. Рассмотрим множество $\EuScript M(\gamma, \eta)$ сфер, которые касаются обеих сфер $\gamma$~и~$\eta$. Заметим что множество $\EuScript M(\gamma, \eta)$ распадается на два класса эквивалентности по типу касаний. Если сфера $\alpha$ касается $\gamma$ и $\eta$ одинаковым образом (обеих внутренним, или обеих внешним), то $\alpha$ принадлежит {\it одному классу}. Если же $\alpha$ касается $\gamma$ и $\eta$ различным образом (одной сферы внутренним, а другой внешним, или наоборот), то $\alpha$ принадлежит {\it другому классу}.
Прямые, проходящие через точки касания $\gamma$ и $\eta$ со сферами одного класса, проходят через общую точку. Для сфер одного класса эта точка~-- один из двух центров инверсии, переводящей $\gamma$ и $\eta$ друг в друга, а для сфер другого класса~-- второй такой центр (эти точки~-- центры подобия сфер $\gamma$ и $\eta$).

\begin{remark}
\label{rem_similitude_center}
Все это имеет место быть и в случае, если, скажем, сфера $\eta$ вырождается в плоскость $\pi$ (сферу бесконечно большого радиуса). Тогда рассмотренные выше инверсные центры $\gamma$ и $\pi$~-- это точки сферы $\gamma$, касательные плоскости в которых параллельны $\pi$.
\end{remark}

Следующая теорема является главным результатом этого параграфа. Она описывает семейство коник $\sigma$, которые вместе с данной окружностью $\Sigma$ образуют {\it 3-пару Понселе} $(\Sigma, \sigma)$, т.е. для них существует треугольник, вписанный в $\Sigma$ и описанный около $\sigma$. Из этой теоремы практически мгновенно следует теорема Грейса, что мы сразу покажем после ее формулировки.

\begin{theorem} [О 3-парах Понселе]
\label{3-pair}
Пусть даны плоскость $\pi$ и окружность $\Sigma$ на ней. Фиксируем сферу $\gamma$, содержащую окружность $\Sigma$, и рассмотрим множество $\EuScript M(\gamma, \pi)$ сфер, касающихся сферы $\gamma$ и~плоскости $\pi$. Тогда если сферы $\alpha$ и $\beta$ пробегают разные классы множества $\EuScript M(\gamma, \pi)$, то описанный около них конус $\EuScript K$ высекает на плоскости $\pi$ семейство коник $\sigma$, образующих 3-пару Понселе с окружностью~$\Sigma$.
\end{theorem}


{\tt Доказательство Теоремы Грейса.}
Пусть $\alpha$ и $\beta$~-- две касательные сферы тетраэдра $ABCD$, гомотетичные относительно вершины $D$. Рассмотрим сферу $\gamma$, касающуюся сфер $\alpha$ и $\beta$ и проходящую через вершины $A$ и $B$. Таких сфер, вообще говоря, целых четыре. Но две из них в данном случае вырождены в плоскости $\langle DAB\rangle$ и $\langle ABC\rangle$, которые принадлежат разным классам множества $\EuScript M(\alpha, \beta)$. Тогда оставшиеся две сферы тоже принадлежат разным классам и в качестве $\gamma$ выберем ту, которая принадлежит другому, нежели плоскость $\langle ABC\rangle$, классу. Пусть она пересекает плоскость $\langle ABC\rangle$ по окружности~$\Sigma$. Описанный около $\alpha$ и $\beta$ конус с вершиной $D$ пересекает плоскость $\langle ABC\rangle$ по конике $\sigma$, касающейся сторон треугольника $ABC$. По Теореме о 3-парах Понселе вершина $C$ также должна лежать на окружности $\Sigma$.

{\hfill $\Box$}
\medskip

{\tt Доказательство Теоремы~\ref{3-pair} о 3-парах Понселе.}


Пусть $F_{\alpha}$ и $F_{\beta}$~-- тоски касания сфер $\alpha$ и $\beta$ с плоскостью $\pi$, которые по теореме Данделена~(1822,~\cite{Dandelin}) являются фокусами коники~$\sigma$. Далее будем считать, что точки $F_{\alpha}$ и $F_{\beta}$ не совпадают друг с другом и с центром окружности $\Sigma$. Эти частные случаи сводятся к общему малым шевелением сфер $\alpha$ и $\beta$ и утверждение теоремы для них получается предельным переходом. Если $I$~-- одна из круговых точек плоскости $\pi$, то $I\in\Sigma$. Обозначим через $P_{\alpha}$ и $P_{\beta}$ точки вторичного пересечения прямых $IF_{\alpha}$ и $IF_{\beta}$ с коникой $\Sigma$. Тогда треугольник $IP_{\alpha}P_{\beta}$ вписан в окружность $\Sigma$, прямые $IP_{\alpha}$ и $IP_{\beta}$ касаются коники $\sigma$, и нам достаточно доказать, в силу теоремы Понселе, что прямая $P_{\alpha}P_{\beta}$ тоже касается коники $\sigma$.

\begin{figure}[htb]
\center
\includegraphics[scale=1]{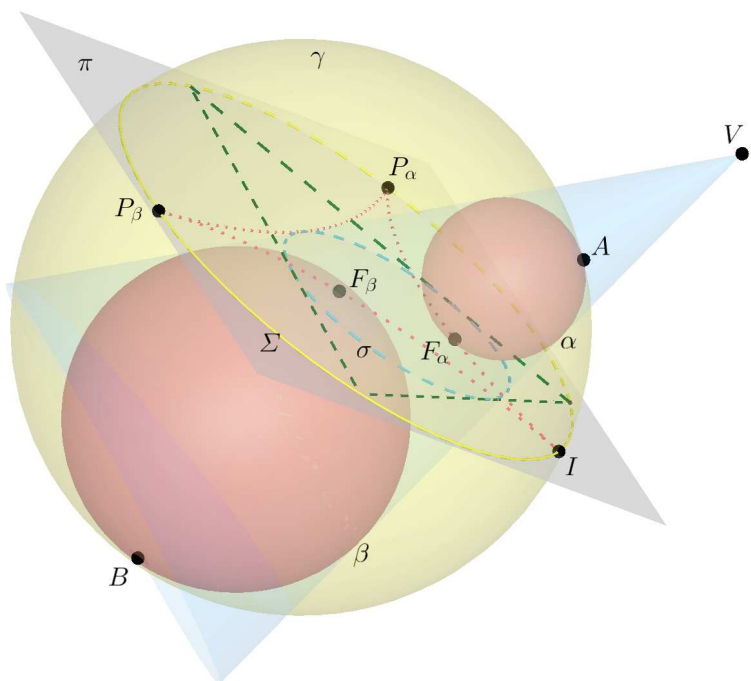}
\caption{\small{\it 3-пары Понселе $(\Sigma, \sigma)$. Мнимые касательные представлены дугообразными розовыми отрезками.}}
\label{fig_3-pairs_Poncelet}
\end{figure}

Пусть $A$ и $B$~-- точки касания сферы $\gamma$ со сферами $\alpha$ и $\beta$. Заметим, что прямая $IF_{\alpha}$ является образующей сферы $\alpha$. Обозначим через $l_A$ одну из двух образующих сферы $\alpha$ в точке $A$, которая пересекает образующую $IF_{\alpha}$ (т.е. $l_A$ и $IF_{\alpha}$ принадлежат разным семействам образующих сферы $\gamma$). Поскольку $l_A$ является также образующей и сферы $\gamma$, точка пересечения $l_A\cap IF_{\alpha}$~-- это одна из двух точек пересечения прямой $IF_{\alpha}$ со сферой $\gamma$, т.е. это либо точка $I$, либо точка $P_{\alpha}$.

Заметим, что первый случай не возможен в силу нашей договоренности считать, что точка $F_{\alpha}$ отлична от центра окружности~$\Sigma$. В самом деле, $I$ лежала бы тогда в пересечении касательных плоскостей сферы $\alpha$ в точках $A$~и~$F_{\alpha}$, т.е. полярно-сопряженная к $AF_{\alpha}$ относительно $\alpha$ прямая содержала бы круговую точку $I$. А так как она вещественная и потому не может быть изотропной, она являлась бы бесконечно-удаленной, т.е. касательные плоскости сферы $\alpha$ в точках $A$ и $F_{\alpha}$ были бы параллельны, а точка $F_{\alpha}$ совпадала бы с центром окружности $\Sigma$.

Таким образом, прямая $AP_{\alpha}$ является общей образующей $l_A$ сфер $\alpha$ и $\gamma$ в точке $A$, и аналогично, прямая $BP_{\beta}$ совпадает с $l_B$~-- одной из двух общих образующих сфер $\beta$ и $\gamma$ в точке $B$. Покажем, что $l_A$ и $l_B$ компланарны.

Для этого рассмотрим гомотетию с центром $A$, переводящую $\alpha$ в $\gamma$. Пусть $g_A$~-- образующая сферы~$\gamma$, в которую переходит образующая $IF_{\alpha}$ сферы $\alpha$. Заметим, что

1) $I\in g_A$, поскольку $g_A\,\|\,IF_{\alpha}$,

2) прямая $g_A$ инцидентна с прямой $l_A$, т.\,к. прямая $l_A$ инвариантна при рассмотренной гомотетии и инцидентна с прямой $IF_{\alpha}$. Т.\,е. $g_A$ и $l_A$~-- две образующие сферы $\gamma$, принадлежащие разным семействам.

Аналогично, если $g_B$~-- образующая сферы~$\gamma$, в которую переходит образующая $IF_{\beta}$ сферы $\beta$ при гомотетии с центром $B$, переводящей $\beta$ в $\gamma$, то

3) $I\in g_B$,

4) $g_B$ и $l_B$~-- тоже две образующие сферы $\gamma$, принадлежащие разным семействам.

Из замечания~\ref{rem_similitude_center} следует, что прямые $g_A$ и $g_B$ проходят через различные инверсные центры сферы $\gamma$ и плоскости $\pi$, а потому различны. Тогда из 1) и 3) следует, что образующие $g_A$ и $g_B$ сферы $\alpha$ имеют общую точку и, значит, принадлежат разным семействам, откуда в силу 2) и 4) следует, что образующие $l_A$ и $l_B$ тоже из разных семейств, а потому компланарны.

Теперь рассмотрим плоскость $\langle l_A; l_B\rangle$, которая в силу утверждения~[iii] Предложения~\ref{prop_generating_lines} касается обеих сфер $\alpha$ и $\beta$. Заметим, что вершина конуса $\EuScript K$ содержит прямую $AB$. Действительно, поскольку конус $\EuScript K$ пересекает $\pi$ по невырожденной конике, его вершина не лежит на $\pi$. Так как $\alpha$ и $\beta$ из разных классов множества $\EuScript M(\gamma, \pi)$, то $\gamma$ и $\pi$ из разных классов множества $\EuScript M(\alpha, \beta)$. Значит, прямая $AB$ проходит через инверсный центр сфер $\alpha$ и $\beta$, который не лежит на плоскости $\pi$.

Т.о., $\langle l_A; l_B\rangle$~-- касательная плоскость конуса $\EuScript K$, а потому пересекает плоскость $\pi$ по прямой, касающейся коники $\sigma$. Осталось заметить, что $\langle l_A; l_B\rangle$ пересекает $\pi$ по прямой $P_{\alpha}P_{\beta}$, и таким образом, треугольник $IP_{\alpha}P_{\beta}$ является вписано-описанным.

{\hfill $\Box$}
\section{Формулы Эйлера-Чаппла и up-in-ex-touch-аналог теоремы Фейербаха}
\label{ch_up-in-ex-touch}
\begin{theorem} [Euler, Chapple]
\label{th_Euler-Chapple}
Пусть $R,\ r$ и $r_a$ -- радиусы описанной, вписанной и вневписанной окружностей произвольного треугольника, $d$ и $d_a$~-- расстояния от центра описанной окружности до центров вписанной и вневписанной. Тогда выполняются следующие соотношения
\begin{equation}\label{eq_Euler-Chapple}
d^2 = R^2 - 2Rr
\end{equation}
\begin{equation}\label{eq_Euler-Chapple2}
d_a^2 = R^2 + 2Rr_a
\end{equation}
\end{theorem}

Мы приведем два, наверное, самых коротких доказательства этой теоремы. Для этого рассмотрим сферу $\Delta$, построенную диаметрально на описанной окружности, наовем ее {\it описанной сферой треугольника}, сферу $\delta$ радиуса $r$, касающуюся плоскости треугольника в центре его вписанной окружности, наовем ее {\it вписано-поднятой}, и сферу $\delta_a$ радиуса $r_a$, касающуюся плоскости треугольника в центре соответствующей вневписанной окружности, наовем ее {\it вневписано-поднятой}.

Заметим, что соотношения~(\ref{eq_Euler-Chapple}),~(\ref{eq_Euler-Chapple2}) можно переписать в виде равенств
$$
d^2 + r^2 = (R-r)^2, \ \ \ d^2 + r_a^2 = (R+r_a)^2,
$$
которые равносильны касанию сфер $\Delta$ и $\delta$, $\Delta$ и $\delta_a$.


\begin{wrapfigure}[15]{r}{190pt}
 \ \ \includegraphics[scale=0.53]{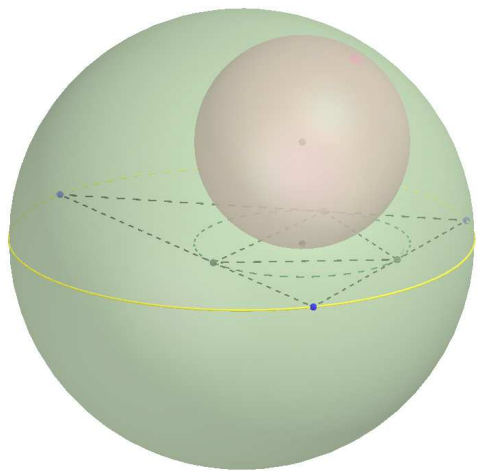}
\caption{\small{\it Сферы $\Delta$ и $\delta$ касаются друг друга}}
\end{wrapfigure}

\noindent\ \ \ \ {\tt Доказательство 1.}
Касания $\Delta$ и $\delta$, $\Delta$ и $\delta_a$ сразу следует из Теоремы Грейса. Действительно, рассмотрим тетраэдр с основанием $ABC$ и вершиной $D$ на бесконечности в перпендикулярном к плоскости $(ABC)$ направлении. Тогда сфера $\delta$ является его вписанной сферой, симметричная ей относительно плоскости $(ABC)$~-- его вневписанной сферой, а следовательно, сфера $\Delta$~-- его сферой Грейса. Для пары $\Delta$ и $\delta_a$ рассуждение аналогично.
{\hfill $\Box$}
\medskip

Это доказательство примечательно своей лаконичностью и красотой, но использование сложной Теоремы Грейса может выглядеть как «стрельба из пушки по воробьям». Поэтому приводим другое

\noindent\ \ \ \ {\tt Доказательство 2.}
Сделаем инверсию относительно сферы, построенной диаметрально на вписанной окружности. Заметим, что сфера $\Delta$ переходит в сферу $\Delta'$, построенную диаметрально на окружности, проходящей через середины сторон треугольника Жергона (вершинами которого являются точки касания вписанной окружности $\triangle ABC$ со сторонами). А сфера $\delta$ переходит в плоскость $\delta'$, удаленную от плоскости $(ABC)$ параллельно на расстояние $\dfrac{\;r}{\;2}$\,. Поскольку, радиус сферы $\Delta'$, очевидно, тоже равен $\dfrac{\;r}{\;2}$, сферы $\Delta'$ и $\delta'$, а следовательно, и сферы $\Delta$ и $\delta$ касаются друг друга.
{\hfill $\Box$}

Заметим, что доказанное свойство касания сферы $\Delta$ с четырьмя сферами $\delta,\, \delta_a,\, \delta_b,\, \delta_c$ является своего рода тоже неким аналогом теоремы Фейербаха в пространстве.
\begin{theorem} [Up-in-ex-touch]
\label{th_up-in-ex-touch}
Описанная сфера треугольника касается его вписано-поднятой и четырех вневписано-поднятых сфер.
\end{theorem}
\begin{figure}[htb]
\center
\includegraphics[scale=0.95]{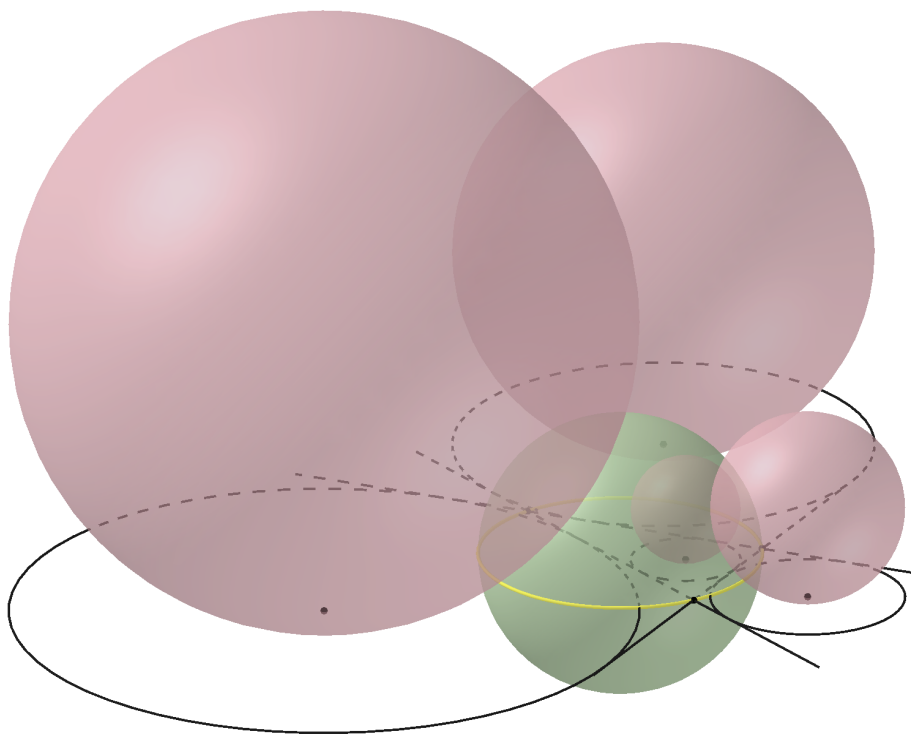}
\caption{\small{\it Up-in-ex-touch-аналог теоремы Фейербаха.}}
\label{fig_Up-in-ex-touch}
\end{figure}

Заметим также, что сфера $\Delta$ касается не только сфер $\delta,\, \delta_a,\, \delta_b,\, \delta_c$, но и еще четырех симметричных им относительно плоскости треугольника, т.е. целых восьми сфер.

\section{Теорема Лагерра и ее применение к тетраэдру}
\label{ch_Laguerre}
\begin{theorem}[Laguerre~\cite{Laguerre}, 1879]
\label{th_Laguerre}
Окружность $\Sigma$ радиуса $R$ с центром в точке $O$ и коника $\sigma$ с фокусами $F_{\alpha},\,F_{\beta}$ и малой полуосью $b$ образуют 3-пару Понселе тогда и только тогда, когда выполняется соотношение
\begin{equation}\label{eq_Laguerre}
(R^2 - d_{\alpha}^2)(R^2 - d_{\beta}^2) = 4R^2b^2,
\end{equation}
где $d_{\alpha} = |OF_{\alpha}|, d_{\beta} = |OF_{\beta}|$.
\end{theorem}
\begin{remark}
Малая полуось $b$ может быть как действительной (у эллипсов), так и мнимой (у гипербол). В первом случае из формулы Лагерра видно, что фокусы эллипса должны лежать либо оба внутри окружности, либо оба вне. Во втором случае, у гиперболы, один фокус должен лежать внутри окружности, другой~-- снаружи.
\end{remark}
\begin{remark}
\label{rem_parabolas}
Если коника $\sigma$ является параболой, то условие существования вписано-описанных треугольников для пары $(\Sigma,\,\sigma)$ становится совсем простым: $d=R$, где $d=|OF|$, т.е. фокус $F$ параболы должен лежать на окружности. Это следует из известной теоремы Ламбера.
\end{remark}

{\tt Доказательство Теоремы Лагерра} ($\Rightarrow$) Пусть $\gamma$ -- произвольная сфера, содержащая окружность~$\Sigma$, а cфера $\alpha$ касается в точке $F_{\alpha}$ плоскости $\pi$, содержащей окружность~$\Sigma$, а~также касается сферы~$\gamma$. Рассмотрим произвольный вписано-описанный треугольник $ABC$ и проведем через его стороны касательные плоскости к сфере $\alpha$. Они пересекаются в некоторой точке $D$, образуя тетраэдр $ABCD$, у которого сфера $\alpha$ является одной из касательных сфер, а $\gamma$~-- сферой Грейса, которая касается также другой касатеьной сферы~$\beta$ тетраэдра $ABCD$, гомотетичной $\alpha$ относительно вершины $D$. Как известно, сферы $\alpha$ и $\beta$ касаются плоскости $\pi$ в точках, изогонально сопряженных относительно $\triangle ABC$. Кроме того, поскольку $F_{\alpha}$ и $F_{\beta}$~-- фокусы вписанной в $\triangle ABC$ коники $\sigma$, они также изогонально сопряжены. Отсюда заключаем, что сфера $\beta$ касается плоскости $\pi$ в точке $F_{\beta}$.

Нам понадобится одна очень простая лемма
\begin{lemma}[Thebault \cite{Thebault}, 1922]
\label{lemma_Thebault}
Для малой полуоси $b$ коники, высекаемой описанным около сфер $\alpha$ и $\beta$ конусом на их общей касательной плоскости, выполняется соотношение
\begin{equation}\label{eq_Thebault}
|b^2| = r_{\alpha} r_{\beta}
\end{equation}
\end{lemma}

Пусть $S_{\alpha}$ и $S_{\beta}$~-- две диаметрально противоположные точки на $\gamma$ в перпендикулярном к плоскости $\pi$ направлении, которые являются инверсными центрами сферы $\gamma$ и плоскости $\pi$ (см. замечание~\ref{rem_similitude_center}). Учитывая, что сферы $\alpha$ и $\beta$ принадлежат разным классам множества $\EuScript M(\gamma, \pi)$ (см. доказательство теоремы Грейса),
легко выразить радиусы сфер $\alpha$ и $\beta$:
\begin{equation}\label{eq_r_alpha_r_beta}
r_{\alpha} = \left|\frac{\Sigma(F_{\alpha})}{2\pi(S_{\alpha})}\right|, \ \ r_{\beta} = \left|\frac{\Sigma(F_{\beta})}{2\pi(S_{\beta})}\right|,
\end{equation}
где $\Sigma(F_{\alpha}) = d_{\alpha}^2 - R^2$ и  $\Sigma(F_{\beta}) = d_{\beta}^2 - R^2$ -- степени точек $F_{\alpha}$ и $F_{\beta}$ относительно окружности $\Sigma$, \linebreak
а~$\pi(S_{\alpha}), \pi(S_{\beta})$~-- расстояния от точек $S_{\alpha}$ и $S_{\beta}$ до плоскости $\pi$.

Перемножим равенства~(\ref{eq_r_alpha_r_beta}) и учтем, что $\pi(S_{\alpha})\pi(S_{\beta}) = R^2$. Получим, что в равенстве~(\ref{eq_Laguerre}) левая и правая части равны по модулю. Правая часть отрицательна только в случае, если коника $\sigma$ является гиперболой. Такое происходит только тогда, когда сферы $\alpha$ и $\beta$ касаются описанного около них конуса с вершиной $D$ по разные стороны от $D$, а плоскости $\pi$~-- по одну сторону. Тогда сферы $\gamma$ они должны касаться по разные стороны, а следовательно, точки $F_{\alpha}$ и $F_{\beta}$ их касания с $\pi$ относительно окружности $\Sigma$ лежат тоже по разные стороны и левая часть~(\ref{eq_Laguerre}) в этом случае также отрицательна. Таким образом, модули можно снять и равенство~(\ref{eq_Laguerre}) считать доказанным.

($\Leftarrow$) Пусть выполняется (\ref{eq_Laguerre}). Если коники $(\Sigma,\,\sigma)$ не образуют 3-пару Понселе, то можно изменить малую полуось $b$ коники $\sigma$ так, чтобы они образовали 3-пару Понселе. Тогда по уже доказанному тоже должно выполняться равенство (\ref{eq_Laguerre}), следовательно величина $b$ не изменилась, т.е. $(\Sigma,\,\sigma)$ как раз и образуют 3-пару Понселе

{\hfill $\Box$}
\medskip

Теорема Лаггера, примененная к тетраэдру, позволяет получить следующее интересное метрическое соотношение для касательных сфер тетраэдра.
\begin{theorem}
\label{th_cos_cos}
Пусть $\Delta_D$~-- сфера, описанная около грани $ABC$ тетраэдра $ABCD$, $\alpha$ и $\beta$~-- две касательные сферы, гомотетичные относительно $D$. Тогда произведение косинусов углов, которые сфера $\Delta_D$ образует с $\alpha$ и $\beta$ (среди них один угол мнимый), равно\  $1$, если $\alpha$ и $\beta$ касаются $\langle ABC\rangle$ с одной стороны, или\ \  $-1$, если с разных. Иными словами,
\begin{equation}\label{eq_cos_cos}
\cos(\widehat{\Delta_D, \alpha})\cos(\widehat{\Delta_D, \beta}) = {\rm sign}\,{k},
\end{equation}
где $k$~-- коэффициент упомянутой гомотетии с центром $D$.
\end{theorem}

{\tt Доказательство} Пусть $O_D$ и $R$~-- центр и радиус сферы $\Delta_D$; \ $r_{\alpha},\,r_{\beta}$~-- радиусы сфер $\alpha$ и $\beta$; \ $\mathfrak{D}_{\alpha}, \mathfrak{D}_{\beta}$~-- расстояния между центрами $\Delta_D$ и $\alpha$, $\Delta_D$ и $\alpha$; $d_{\alpha}, d_{\beta}$~-- расстояния от $O_D$ до точек $F_{\alpha}$ и $F_{\beta}$ касания плоскости $\langle ABC\rangle$ со сферами $\alpha$ и $\beta$. Пусть $\Sigma = \odot(ABC)$, а конус
$\EuScript K$ с вершиной $D$, описанный около $\alpha$ и $\beta$, пересекает плоскость $\langle ABC\rangle$ по конике $\sigma$.

Воспользуемся леммой~\ref{lemma_Thebault} и заметим, что в нашей конструкции с тетраэдром равенство (\ref{eq_Thebault}) можно уточнить
\begin{equation}\label{eq_Thebault_sign}
b^2 = r_{\alpha} r_{\beta}\,{\rm sign}\,{k},
\end{equation}
поскольку $b^2$ может быть отрицательным, только если коника $\sigma$ является гиперболой, что возможно лишь в том случае, если вершина конуса $\EuScript K$ является центром отрицательной гомотетии сфер $\alpha$ и $\beta$, т.е. они вписаны в $\EuScript K$ по разные стороны от его вершины.

По теореме Лагерра для пары $(\Sigma, \sigma)$ имеем
\begin{equation}\label{eq_Laguerre_tetrahedron}
(R^2 - d_{\alpha}^2)(R^2 - d_{\beta}^2) = 4R^2b^2.
\end{equation}
По теореме Пифагора
$$d_{\alpha}^2 = \mathfrak{D}_{\alpha}^2 - r_{\alpha}^2,\ \ d_{\beta}^2 = \mathfrak{D}_{\beta}^2 - r_{\beta}^2$$
Подставляя эти равенства и~(\ref{eq_Thebault_sign}) в соотношение (\ref{eq_Laguerre_tetrahedron}), получаем требуемое соотношение

$$\frac{R^2 + r_{\alpha}^2 -  \mathfrak{D}_{\alpha}^2}{2 R\, r_{\alpha}}\cdot \frac{R^2 + r_{\beta}^2 -  \mathfrak{D}_{\beta}^2}{2 R\, r_{\beta}} = {\rm sign}\,{k}$$
{\hfill $\Box$}
\medskip

\section{Трехмерный аналог формулы Эйлера-Чаппла}
\label{ch_R3_Euler-Chapple}
В связи с теоремой Эйлера-Чаппла возникает естественный вопрос о возможности ее трехмерного обобщения на случай тетраэдра.
Этот вопрос был поставлен впервые Ж.\,Д.\,Жергонном в 1816 году в издаваемом им журнале\footnote{{\it {Annales de math\'{e}matiques pures et appliqu\'{e}es}},\,\ {\bf 6} (1815-1816), p. 228.} в виде краткой сноски, относящейся к тетраэдру с радиусами описанной сферы $R$, вписанной~-- $r$ и расстоянием $d$ между их центрами:

\begin{center}
\begin{minipage}[h]{0.85\linewidth}
{\it «Il serait sur-tout int\'{e}ressant de savoir si \,$d$\, peut \^{e}tre exprim\'{e} uniquement \linebreak en fonction de $R$ et $r$.
\hfill {\it J.\,D.\,G.}»}
\end{minipage}
\end{center}

Спустя восемь лет в том же журнале было опубликовано положительное решение этой задачи в работе Дюрранда~\cite{Durrande}, где он доказал следующее соотношение:
\begin{equation}\label{eq_Durrande}
d^2 = (R + r)(R - 3r).
\end{equation}

Этот результат получил широкое признание и в течение многих лет на него ссылались в литературе, например, в таких авторитетных изданиях как Математическая энциклопедия Клейна «Encyklop\"{a}die der mathematischen Wissenschaften»~\cite{Zacharias} (первая в мире математическая энциклопедия) и «Enciclopedia delle matematiche elementari»~\cite{Biggiogero} (крупнейшая энциклопедия по математике, изданная в Италии). Однако, формула Дюрранда~(\ref{eq_Durrande}) оказалась неверной, а ответ на вопрос Жергонна~-- отрицательным: не существует общей для всех тетраэдров функциональной зависимости между $R, r$ и $d$. Доказательство Дюрранда было практически безупречным, но незаметная ошибка заключалась в его убежденности, что описанная и вписанная сфера непременно должны иметь такую зависимость. Вопрос Жергонна можно было бы сформулировать так: {\it каковы условия существования вписано-описанного тетраэдра для двух данных сфер?}

Оказывается никаких необходимых условий для этого не требуется.

\begin{theorem}
\label{th_circum-inscribes_tetra}
Для любых двух невырожденных квадрик общего положения существует бесконечное семейство вписано-описанных тетраэдров. Любая касательная плоскость ко вписанной квадрике может содержать грань такого тетраэдра, а его вершиной может быть произвольная точка описанной квадрики.
\end{theorem}

В работе Фонтене~\cite{Fontene} 1899 года эта теорема считается уже известной (см. также~\cite{Grace}).

Итак, в отличие от плоского случая, в пространстве для любых двух произвольных сфер всегда существует вписано-описанный в них тетраэдр, причем он может динамически вращаться около этих сфер, все время оставаясь вписано-описанным. При этом, любая точка описанной сферы может быть вершиной такого тетраэдра.

Но не для любых двух сфер такой тетраэдр может быть вещественным. Он всегда существует, но может быть мнимым. Критерием существования {\it вещественного} вписано-описанного тетраэдра является следующее условие Грейса, исправляющее соотношение Дюрранда~(\ref{eq_Durrande}):

\begin{theorem}[Grace~\cite{Grace}, 1917]
\label{th_Grace_inequality}
Для данных двух сфер $S$ и $T$ критерием существования вещественного вписано-описанного тетраэдра, у которого вершины лежат на $S$, а плоскости граней касаются $T$, являются следующие условия для каждого случая взаимного расположения $S$ и $T$:
\begin{description}
 \item[(a)] $T$ вложена в $S$ и $$d^2 \leqslant (R + r)(R - 3r);$$
 \item[(b)] $T$ и $S$ расположены одна вне другой;
 \item[(c)] $T$ и $S$ пересекаются по действительной окружности и
  $$d^2 \leqslant (R - r)(R + 3r).$$
\end{description}
\end{theorem}


\section{Вращение Понселе вписано-описанного тетраэдра}
\label{ch_R3_Ponc}
Теорема~\ref{th_circum-inscribes_tetra} позволяет рассмотреть динамику «вращения» вписано-описанного тетраэдра. Эта динамика не столь однозначна, как в плоской теореме Понселе. Это показывает следующая теорема.
\begin{theorem}[\cite{Grace}]
\label{th_dynamic_tetra}
Пусть вершины тетраэдра лежат на квадрике $S$, а грани касаются квадрики $T$. Тогда при фиксации плоскости $\pi$ одной из его граней противоположная вершина $P$ может при этом варьироваться, пробегая сечение квадрики $S$ некоторой плоскостью $\pi'$.
\end{theorem}

Таким образом, тетраэдр вращается с намного большей свободой, чем вписано-описанный многоугольник. Когда выбрана плоскость $\pi$, существует целая коника для выбора произвольной точки на ней в качестве вершины $P$, а для каждой такой пары $P$ и $\pi$ существует однопараметрическое семейство вписано-описанных треугольников, каждый из которых может быть противоположной к вершине $P$ гранью вписано-описанного тетраэдра. Таким образом, в общем случае существует $4$-параметрическое семейство тетраэдров.

У плоской теоремы Понселе есть такой «эффект замыкания»: если начиная с некоторой начальной точки $A_1$ строится последовательно вписано-описанная ломаная $A_1A_2\ldots A_n$ и оказывается, что звено $A_1A_n$ тоже касается вписанной коники, замыкая ломаную, то такое замыкание будет происходить всегда.

Если же по аналогии строить вписано-описанный тетраэдр для двух данных квадрик $S$ и $T$, последовательно выбирая касательные плоскости его граней, то возникает следующий вопрос. {\it Когда мы провели уже три плоскости, которые образовали вписано-описанный трехгранный угол, всегда ли можно его замкнуть четвертой плоскостью, чтобы образовался вписано-описанный тетраэдр?} Очевидно, что четвертая грань однозначно определяется по первым трем. Критерий того, что она будет тоже касаться вписанной сферы, дает следующая теорема Фонтене.
\begin{theorem}[Fonten\'{e}~\cite{Fontene}]
\label{th_Fontene}
Последовательный процесс построения вписано-описанного тетраэдра замыкается всегда, если и только если квадрики $S$ и $T$ имеют четыре общие образующие. \end{theorem}

Плоскость $\pi$ и вершина $P$ в этом случае могут быть выбраны совсем произвольно и, таким образом, существует 5-параметрическое семейство вписано-описанных тетраэдров. 
\begin{theorem}
\label{th_dynamic_Grace}
Пусть фиксированы описанная сфера $S$ тетраэдра и одна из восьми его касательных сфер $T$, а тетраэдр динамически «вращается» около них, оставаясь вписано-описанным. Тогда все четыре касающиеся $T$ сферы Грейса все время касаются некоторой фиксированной сферы, концентричной с описанной сферой S.
\end{theorem}

{\tt Доказательство.}
Пусть сфера Грейса $G$ проходит через вершины грани $a$ и пусть вписанная сфера $S$ касается сферы $G$ в точке $P$, а плоскости $\langle a\rangle$~-- в точке $Q$. Обозначим центры сфер $S$ и $T$ через $O_S$ и $O_T$. Прямая $PQ$ при вращении тетраэдра проходит через фиксированную точку~-- предельную точку $K$ пучка сфер $\langle S,\, T\rangle$. Кроме того, на прямой $PQ$ лежит инверсный центр $E$ сферы $G$ и плоскости $\langle a\rangle$, касательная в котором к $G$ параллельна плоскости $\langle a\rangle$. Следовательно, $O_SE\|O_TQ$ и $\triangle O_SEK\sim \triangle O_TQK$, откуда получаем такое выражение
$$O_SE = \frac{O_SK}{O_TK}\cdot r_T,$$
правая часть которого является величиной постоянной при вращении тетраэдра. Тогда, сфера с радиусом, равным этой величине, и центром в точке $O_S$ касается сферы Грейса в любой момент вращения.

{\hfill $\Box$}
\medskip
\section{Доказательство теоремы Фейербаха через выход в пространство}
\label{ch_Feuer_into_R3}
\ \ \ \ Пусть $\delta$~-- вписанная окружность треугольника $ABC$ с центром в точке $I$ и радиусом $r$, $H$~-- ортоцентр треугольника $ABC$, точки $A_1, B_1, C_1, I_1$~-- середины отрезков $AH, BH, CH, IH$ ($I_1$~-- инцентр $\triangle A_1B_1C_1$). Описанная около $\triangle A_1B_1C_1$ окружность $\vartheta$~-- это окружность девяти точек $\triangle ABC$.
Пусть также $\odot_{a},\,\odot_{b},\,\odot_{c}$~-- окружности с диаметрами $BC, CA, AB$, $\Delta$ и $\Theta$~-- сферы, построенные диаметрально на окружностях $\delta$ и $\theta$.

{\tt Доказательство Теоремы Фейербаха.}

Заметим, что касание окружностей $\delta$ и $\theta$ равносильно касанию сфер $\Delta$ и $\Theta$.
По Теореме~\ref{th_up-in-ex-touch} для $\triangle A_1B_1C_1$ его описанная сфера $\Theta$ касается его вписано-поднятой сферы $\Upsilon$.
Поэтому касание $\Theta$ и $\Delta$ равносильно тому, что сфера $\Theta$ инвариантна при инверсии, переводящей сферы $\Delta$ и $\Upsilon$ друг в друга.
Заметим, что центр $S$ этой инверсии расположен над точкой $H$ на высоте~$r$ (т.е. $SH\bot(ABC),\, |SH| = r$), а коэффициент инверсии (квадрат радиуса сферы инверсии) равен $|IH|\cdot|I_1H| = \dfrac{\,|IH|^2}{2}$.
Таким образом, достаточно доказать равенство $\Theta(S) = \dfrac{|IH|^2}{2}$, которое в силу того, что $\Theta(S) = \theta(H) + r^2$, равносильно соотношению
\begin{equation}\label{eq_IH}
|IH|^2 - 2r^2 = 2\theta(H)
\end{equation}

Заметим, что левая часть равенства~(\ref{eq_IH}) равна степени точки $H$ относительно окружности $\xi$ радиуса $r\sqrt2$ с центром $I$ ($\xi$ высекает на сторонах $\triangle ABC$ равные отрезки длины~$2r$). Осталось воспользоваться следующим замечательным свойством окружности $\xi$.
\begin{theorem}
\label{th_xi(a)(b)(c)}
Окружности $\xi, \,\odot_{a},\,\odot_{b},\,\odot_{c}$ имеют общий радикальный центр в точке~$H$.
\end{theorem}
Тогда заметим, что степень точки $H$ относительно окружности $\theta$ в два раза меньше ее степени относительно окружностей $\odot_{a},\,\odot_{b},\,\odot_{c}$ и равенство~(\ref{eq_IH}) равносильно утверждению $\xi(H) = \odot_{a}(H) = \odot_{b}(H) = \odot_{c} (H)$ Теоремы~\ref{th_xi(a)(b)(c)}.

{\hfill $\Box$}
\medskip

Для доказательства Теоремы~\ref{th_xi(a)(b)(c)} рассмотрим окружность $\chi_a$, диаметром которой является жергониана вершины $A$ (т.е. отрезок, соединяющий $A$ с точкой касания вписанной окружности $\delta$ со стороной $BC$) и воспользуемся следующим свойством окружности $\chi_a$, возможно, имеющим и самостоятельный интерес.
\begin{lemma}[$\chi_a$-лемма]
\label{chi_a}
Окружности $\chi_a,\, \xi,\, \odot_{a}$ принадлежат одному пучку.
\end{lemma}

{\tt Доказательство Теоремы~\ref{th_xi(a)(b)(c)}.}
Достаточно проверить, что $H\in\rm{rad}(\xi, \odot_{a})$. \\
Заметим, что $\rm{rad}(\chi_a, \odot_{b})$~-- это высота $AH$, $\rm{rad}(\odot_{a}, \odot_{b})$~-- это высота $CH$, следовательно,\\
$$H = \rm{rad}(\chi_a, \odot_{a}, \odot_{b})\in \rm{rad}(\chi_a, \odot_{a}) = \rm{rad}(\xi, \odot_{a}),$$
где последнее равенство верно в силу $\chi_a$-леммы.
{\hfill $\Box$}
\medskip

{\tt Доказательство $\chi_a$-леммы.}

Воспользуемся следующим известным метрическим соотношением для пучков окружностей.
\begin{lemma}[О пучке]
\label{lemma_pencil}
Если окружности $\alpha, \beta, \gamma$ лежат в одном пучке, то для любой точки $P\in\gamma$ отношение ее степеней относительно $\alpha$ и $\beta$ постоянно, причем
\begin{equation}
\label{eq_pencil}
\frac{\alpha(P)}{\beta(P)}=\frac{d_{\alpha\gamma}}{d_{\beta\gamma}},
\end{equation}
где $d_{\alpha\gamma}$ и $d_{\beta\gamma}$~-- расстояния между центрами $\alpha, \gamma$ и $\beta, \gamma$.
\end{lemma}

Верно и обратное утверждение.
\begin{lemma}[Обратная лемма о пучке]
\label{lemma_pencil_converse}
 Пусть центры окружностей $\alpha, \beta, \gamma$ коллинеарны, и на окружности $\gamma$ имеется такая точка $P$, для которой выполняется соотношение~(\ref{eq_pencil}). Тогда окружности $\alpha, \beta, \gamma$ принадлежат одному пучку.
\end{lemma}


\begin{figure}[htb]
\center
\includegraphics[scale=0.6]{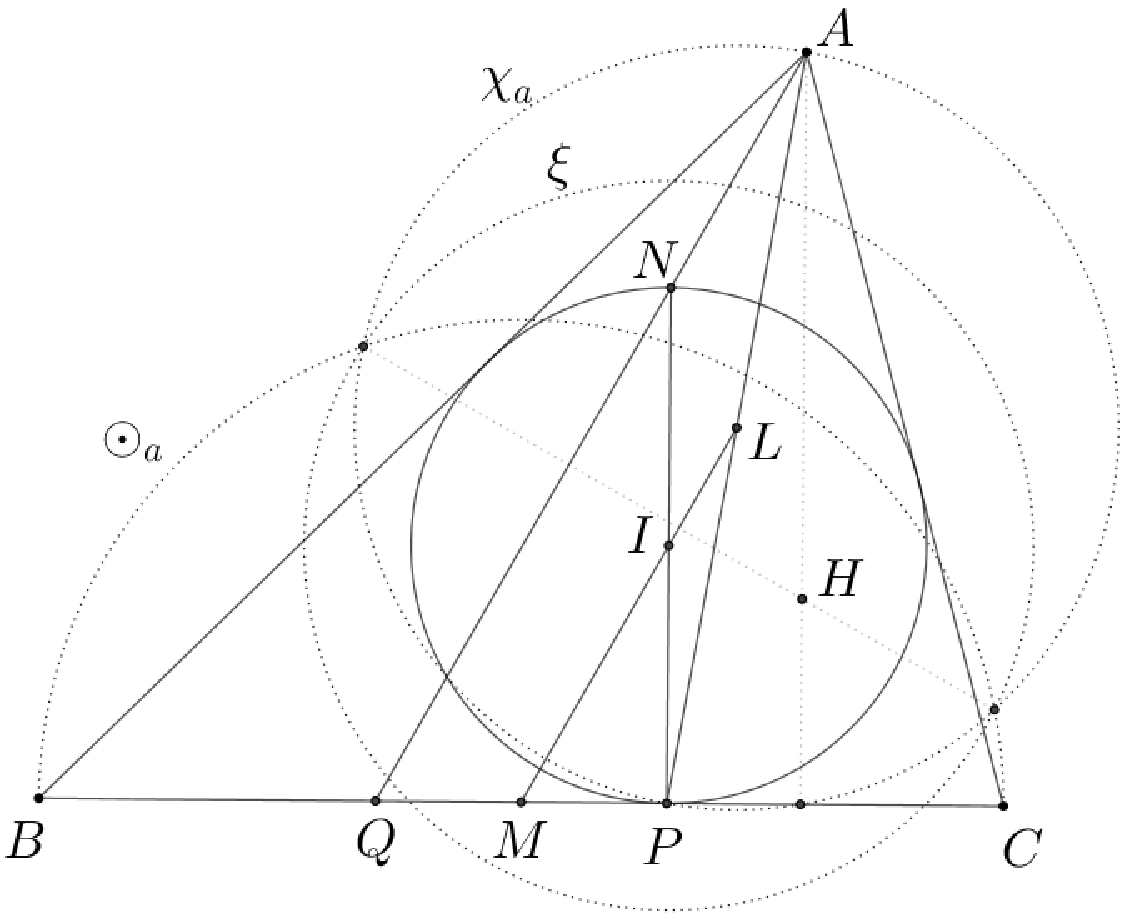}
\caption{\small{\it Окружности $\xi, \chi_a, \odot_{a}$ принадлежат одному пучку}}
\label{fig_chi_a_lemma}
\end{figure}

В качестве окружностей $\alpha, \beta, \gamma$ из Обратной леммы о пучке возьмем окружности $\odot_{a},\, \xi,\, \chi_a$, центры $M, I, L$ которых лежат на средней линии $ML$ треугольника $APQ$. При этом,
$$\frac{LM}{LI} = \frac{AQ}{AN} = \frac{r_a}{r}.$$
Для точки $P\in\chi_a$ имеем $\alpha(P) = - (p-b)(p-c),\ \ \beta(P) = - r^2.$
Тогда~(\ref{eq_pencil}) запишется в виде соотношения
$$\frac{(p-b)(p-c)}{r^2} = \frac{r_a}{r},$$
которое равносильно легко проверяемому равенству
$$(p-b)(p-c) = r\,r_a.$$

{\hfill $\Box$}
\medskip

{\tt Доказательство $\chi_a$-леммы выходом в пространство.}
Заметим, что окружности $\chi_a,\ \delta$ и окружность $\odot_{PQ}$ с диаметром на отрезке $PQ$ лежат в одном пучке. Поднимем их центры перпендикулярно плоскости $\langle ABC\rangle$, сохраняя коллинеарность: $L\to\overline{L}, I\to\overline{I}, M\to\overline{M}$, и пусть $\overline{L}L = \dfrac{r}{2}, \overline{I}I = r$. Тогда легко найти, что $\overline{M}M = \dfrac{r + r_a}{2}$. При этом сферы $S(\overline{L}), S(\overline{J}), S(\overline{M})$ с центрами $\overline{L}, \overline{I}, \overline{M}$, содержащие окружности $\chi_a,\ \delta,\ \odot_{PQ}$ соответственно, также принадлежат одному пучку. Рассмотрим плоскость $\pi\|\langle ABC\rangle$, проходящую через $\overline{I}$, и ортогональную проекцию $\triangle ABC\to\triangle A'B'C'$ на плоскость $\pi$. Осталось заметить, что сечениями сфер $S(\overline{L}), S(\overline{J}), S(\overline{M})$ плоскостью $\pi$ являются окружности $\xi', \chi_a', \odot_{a}'$. Действительно, для сечений $S(\overline{L}), S(\overline{I})$ это очевидно, а для $S(\overline{M})$ это легко проверить, поскольку квадрат радиуса окружности ее сечения плоскостью $\pi$ равен
$$|MP|^2 + \left(\frac{r_a + r}{2}\right)^2 - \left(\frac{r_a - r}{2}\right)^2 = |MP|^2 + r_ar = |MP|^2 + (p-b)(p-c) = |MP|^2 + |BP|\cdot|CP| = \left(\frac{\,a}{\,2}\right)^2.$$

Так как при пересечении сфер пучка плоскостью получается пучок окружностей, то $\xi', \chi_a', \odot_{a}'$ принадлежат одному пучку.
{\hfill $\Box$}
\medskip

\section*{Благодарности.}
Автор благодарен В.\,Ю. Протасову за внимание к работе, полезные советы и помощь в подготовке текста.
{\small{}}
\end{document}